\documentclass[
global_equation_numbers,
global_theorem_numbers,
global_figure_numbers,
global_table_numbers,
framed_theorems,
framed_proofs,
print,
]{OCPreprint}

\usepackage{enumitem}
\usepackage{pgfplotstable}

\newcommand\sign{\operatorname{sign}}

\titlehead{}
\title{Maximal monotone operators with non-maximal graphical limit}
\author{Gerd Wachsmuth\footnote{%
		Brandenburgische Technische Universität Cottbus--Senftenberg,
		Institute of Mathematics,
		03046 Cottbus,
		Germany,
		\email{wachsmuth@b-tu.de},
		\url{https://www.b-tu.de/fg-optimale-steuerung/team/prof-gerd-wachsmuth},
		ORCID: 0000-0002-3098-1503%
	}
}
\publishers{\vspace*{-\baselineskip}}
{
	\makeatletter
	\def\and{ and }
	\def\footnote#1{}
	\hypersetup{
		pdftitle={\@title},
		pdfauthor={\@author}
	}
	\makeatother
}

\begin{document}
\maketitle
\begin{abstract}
	We present a counterexample
	showing that the graphical limit of maximally monotone operators
	might not be maximally monotone.
	We also characterize the directional differentiability
	of the resolvent of an operator $B$
	in terms of existence and maximal monotonicity
	of the proto-derivative of $B$.
\end{abstract}

\begin{keywords}
	Maximal monotone operator,
	graphical limit,
	proto-derivative,
	directional differentiability
\end{keywords}

\begin{msc}
	\mscLink{49J53},
	\mscLink{49J52},
	\mscLink{49K40},
	\mscLink{47H04}
\end{msc}

\section{Introduction}
By means of a counterexample
we show that
the (non-empty) graphical limit of maximally monotone operators
may fail to be maximally monotone.
This was raised as an open question in
\cite[Remark~6]{AdlyRockafellar2020}.
We also shed some light on the relation of proto-differentiability
of an operator
and
directional differentiability of its resolvent.
Throughout this work,
we use standard notation,
see, e.g., \cite{AdlyRockafellar2020}.

\section{Graphical limits of maximally monotone operators}
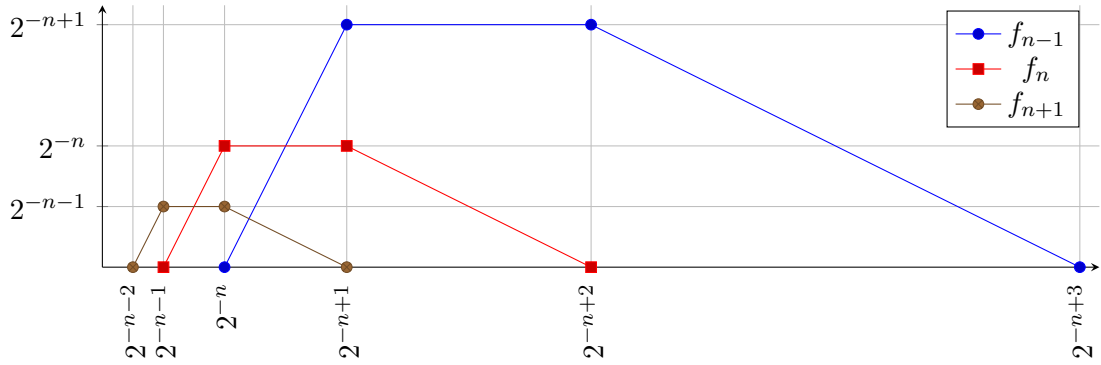
\begin{figure}[ht]
	\centering
	\begin{tikzpicture}
		\begin{axis}[grid=both,
				xmin=0.0,
				xmax=0.51,
				ymax=0.135,
				axis lines=middle,
				width=1.0\textwidth,
				axis equal image,
				xtick={
					{2^-6},
					{2^-5},
					{2^-4},
					{2^-3},
					{2^-2},
					{2^-1}
				},
				xticklabels={
					$2^{-n-2}$,
					$2^{-n-1}$,
					$2^{-n}$,
					$2^{-n+1}$,
					$2^{-n+2}$,
					$2^{-n+3}$
				},
				xticklabel style={rotate=90},
				ytick={
					{2^-5},
					{2^-4},
					{2^-3}
				},
				yticklabels={
					$2^{-n-1}$,
					$2^{-n}$,
					$2^{-n+1}$
				},
			]
			\pgfplotstableread[]{
				0.062500000000000                   0
				0.125000000000000   0.125000000000000
				0.250000000000000   0.125000000000000
				0.500000000000000                   0
			}\fn
			\pgfplotstableread[]{
				0.031250000000000                   0
				0.062500000000000   0.062500000000000
				0.125000000000000   0.062500000000000
				0.250000000000000                   0
			}\fnn
			\pgfplotstableread[]{
				0.015625000000000                   0
				0.031250000000000   0.031250000000000
				0.062500000000000   0.031250000000000
				0.125000000000000                   0
			}\fnnn
			\addplot table{\fn};
			\addplot table{\fnn};
			\addplot table{\fnnn};
			\legend{$f_{n-1}$,$f_n$,$f_{n+1}$}
		\end{axis}
	\end{tikzpicture}
	\caption{Plot of the functions $f_n$.}
	\label{fig:fn}
\end{figure}

For all $n \in \N$, we define the auxiliary function
$f_n \colon [0,\infty) \to \R$
via
\begin{equation*}
	f_n(t)
	:=
	\begin{cases}
		0 & \text{if } t \le 2^{-n-1} \text{ or } t \ge 2^{-n+2}, \\
		2 \parens*{t - 2^{-n-1}} & \text{if } 2^{-n-1} \le t \le 2^{-n}, \\
		2^{-n} & \text{if } 2^{-n} \le t \le 2^{-n+1}, \\
		\frac12 \parens*{2^{-n+2} - t} & \text{if } 2^{-n+1} \le t \le 2^{-n+2}. \\
	\end{cases}
\end{equation*}
Each $f_n$ is globally Lipschitz continuous with Lipschitz constant $2$,
see also
\cref{fig:fn}.

Let $\ell^2$ be the Hilbert space
of square-summable sequences
(an analogue construction is possible in every infinite-dimensional Hilbert space).
We define $T \colon \ell^2 \to \ell^2$
via
\begin{equation*}
	T(x)
	:=
	\sum_{n = 1}^\infty f_n(\norm{x}) e_n,
\end{equation*}
where $\seq{e_n}_n$ is the canonical orthonormal basis of $\ell^2$.
This operator is well defined,
since for each $x \in \ell^2$,
the sum contains at most three non-vanishing terms.
\begin{lemma}
	\label{lem:T_Lipschitz}
	The operator $T$ is globally Lipschitz continuous on $\ell^2$
	with Lipschitz constant at most $\sqrt{17}/2$.
	Moreover,
	for each $x \in \ell^2$
	with $\norm{x} \le 1$
	we have $\norm{T(x)} \ge \norm{x}/2$.
\end{lemma}
\begin{proof}
	First of all, it can be checked easily that $T$
	is continuous on $\ell^2$.
	Now,
	let $x,y \in \ell^2$ be given
	and we denote
	$(x,y) := \set{\lambda x + (1-\lambda) y \given \lambda \in (0,1)}$.

	We are going
	to utilize the mean value inequality from
	\cite[Theorem 2.7]{Penot2013}.
	Since the functions $f_n$ are directionally differentiable,
	$T$ is directionally differentiable on $(x,y) \setminus \set{0}$
	and
	\begin{equation*}
		T'(z; y - x) = \sum_{n = 1}^\infty f_n'( \norm{z}; \dual{z}{y - x} / \norm{z} ) e_n
		\qquad\forall z \in (x,y) \setminus \set{0}.
	\end{equation*}
	This yields the estimate
	\begin{equation*}
		\norm{T'(z; y - x)}^2
		\le
		\norm{y - x}^2
		\sum_{n = 1}^\infty \abs{f_n'(\norm{z}; \sign \dual{z}{y-x})}^2
		\qquad\forall z \in (x,y) \setminus \set{0}.
	\end{equation*}
	By construction of $f_n$,
	the sum contains at most two distinct addends from $\set{\frac12, 2}$.
	Thus,
	\begin{equation*}
		\norm{T'(z; y - x)} \le \sqrt{2^2 + 1/2^2} \norm{y - x} = \sqrt{17}/2 \, \norm{y - x}
		\qquad\forall z \in (x,y) \setminus \set{0}.
	\end{equation*}
	Now,
	\cite[Theorem~2.7]{Penot2013}
	(together with the remark afterwards)
	yields the estimate
	$\norm{T(y) - T(x)} \le \frac{\sqrt{17}}{2} \norm{y - x}$
	for all $x,y \in \ell^2$.

	For every $t \in [0,1]$, there exists $n \in \N$
	such that $f_n(t) \ge t/2$, cf.\ \cref{fig:fn}.
	This implies the second claim.
\end{proof}

Combining \cref{lem:T_Lipschitz}
with
\cite[Example~20.26]{BauschkeCombettes2011},
we find that
$\Id + \alpha T$
is maximally monotone
for all $\alpha \in \R$ with $\abs{\alpha} \le 2/\sqrt{17}$.
For an arbitrary $\alpha$ in this range,
we set $B := \Id + \alpha T$.
For all $m \in \N$,
we define the operator
$B_m \colon \ell^2 \to \ell^2$ via
\begin{equation*}
	B_m(x)
	:=
	m B(x / m)
	=
	x + \alpha m T(x / m).
\end{equation*}
It is easy to check that all the operators $B_m$
are again maximally monotone.
However,
their graphical limit fails to be maximally monotone
in an extreme way.
\begin{theorem}
	\label{thm:graphical_limit}
	Let the maximally monotone operators $B_m \colon \ell^2 \to \ell^2$ be given as above.
	Then,
	the graphical limit of $B_m$ as $m \to \infty$ is the operator
	$Z \colon \ell^2 \mto \ell^2$,
	defined via
	$\graph(Z) = \set{(0,0)}$.
\end{theorem}
\begin{proof}
	We start by the computation of the outer limit of $\graph(B_m)$.
	For $(x,y) \in \limsup_{m \to \infty} \graph(B_m)$,
	we find a sequence $\seq{(x_{m_k},y_{m_k})}_k$ with $(x_{m_k}, y_{m_k}) \in \graph(B_{m_k})$
	and $x_{m_k} \to x$, $y_{m_k} \to y$.
	In particular, we have
	\begin{equation*}
		y_{m_k} = x_{m_k} + \alpha m_k T(x_{m_k} / m_k)
		.
	\end{equation*}
	Since $\seq{x_{m_k}}_k$ is bounded, we have $x_{m_k} / m_k \to 0$.
	Now,
	the structure of $T$ implies that
	\begin{equation*}
		\bracks{T(x_{m_k} / m_k)}_n = 0
		\qquad\text{for $k$ large enough}
	\end{equation*}
	for each fixed $n$.
	Since $y_{m_k} - x_{m_k} = \alpha m_k T(x_{m_k} / m_k)$
	converges,
	the limit can only attain the value $0$
	and, thus, we have $x = y$.
	Since $x_{m_k} / m_k \to 0$,
	we know $\norm{m_k T(x_{m_k} / m_k) } \ge \norm{x_{m_k}}/2$.
	Together with $y_{m_k} - x_{m_k} = \alpha m_k T(x_{m_k} / m_k) \to 0$,
	this gives $x_{m_k} \to 0$.
	Thus, $(x,y) = (0,0)$ is the only point in
	$\limsup_{m \to \infty} \graph(B_m)$.
	Moreover, $(0,0) \in \graph(B_m)$
	shows
	that the limit of $\graph(B_m)$ is $\set{(0,0)}$.
\end{proof}
Clearly, the same argument can be used
for
the operators $B_\tau \colon \ell^2 \to \ell^2$, $\tau \in (0,1)$,
defined via
$B_\tau(x) = \tau^{-1} B(\tau x)$
and for the limiting process $\tau \searrow 0$.
Note that $B_\tau$ is just the finite difference
appearing in the definition of the proto-derivative
of $B$ at $0$ relative to $0$.

\begin{corollary}
	\label{cor:proto}
	The maximal monotone mapping $B$ is proto-differentiable
	at $0$
	and the proto-derivative at $0$ relative to $0 = B(0)$
	is given by the non-maximally monotone operator $Z$
	from \cref{thm:graphical_limit}.
\end{corollary}

\section{Directional differentiability of resolvents}
Let $H$ be a (real) Hilbert space.
For a maximally monotone $B \colon H \mto H$,
we denote by $J_B \colon H \to H$
its single-valued resolvent,
i.e., $J_B := (\Id + B)^{-1}$.
The next result characterizes
the directional differentiability of $J_B$.
\begin{theorem}
	\label{thm:charact}
	Let $B \colon H \mto H$ be maximally monotone.
	For $y \in H$ set $x := J_B(y)$.
	Then, the following are equivalent.
	\begin{enumerate}[label=(\roman*)]
		\item
			$B$ is proto-differentiable at $x$
			relative to $y - x \in B(x)$
			and the proto-derivative $D_p B(x \mathbin{|} y-x) \colon H \mto H$ is maximally monotone,
		\item
			$J_B$ is directionally differentiable at $y$,
			i.e.,
			the limit
			$
				J_B'(y; h) = \lim_{\tau \searrow 0}\frac{J_B( y + \tau h ) - J_B(y)}{\tau}
				$
			exists for all $h \in H$.
	\end{enumerate}
\end{theorem}
\begin{proof}
	``$\Rightarrow$'' follows from \cite[Theorem~1]{AdlyRockafellar2020}
	by setting $A(t,x) := x$, $B(t,x) := B(x)$, $\xi(t) := y + t \, h$.

	``$\Leftarrow$'':
	From \cite[Remark~5]{AdlyRockafellar2020},
	we get that $J_B$ is proto-differentiable at $y$
	for $x = J_B(y)$.
	Consequently,
	\cite[Lemma~2]{AdlyRockafellar2020}
	implies that $B$ is proto-differentiable at $x$
	relative to $y - x$.
	Moreover,
	we get
	the formula
	\begin{equation*}
		D_p J_B(y \mathbin{|} x)
		=
		\set{
			J'(y; \cdot)
		}
		=
		\parens*{
			\Id
			+
			D_p B(x \mathbin{|} y - x)
		}^{-1}
	\end{equation*}
	linking the
	derivatives of $B$ and $J_B$.
	This shows that the resolvent of
	the monotone operator
	$D_p B(x \mathbin{|} y - x)$
	is single-valued.
	By Minty's theorem \cite[Theorem~21.1]{BauschkeCombettes2011},
	we find that $D_p B(x \mathbin{|} y - x)$
	is maximally monotone.
\end{proof}
\cref{thm:graphical_limit,thm:charact} show
that the requirement of the proto-derivative of $B$
being maximally monotone in \cite[Theorem~1]{AdlyRockafellar2020}
cannot be dropped.
The same result can be proved in the $t$-dependent case
considered in \cite{AdlyRockafellar2020}.

\textbf{Acknowledgement.}
{\small%
This work is supported by the DFG Grant
WA 3636/4-2 within the Priority Program 1962 (Non-smooth
and Complementarity-based Distributed Parameter Systems: Simulation and Hierarchical Optimization).
}

\renewcommand*{\bibfont}{\small}
\printbibliography

\end{document}